\documentclass[12pt]{article}
\usepackage{amssymb}
\usepackage{amsmath}
\newcommand{\qedc}{\hfill {$\Box$}}
\newtheorem{thm}{Theorem}[section]
\newtheorem{lem}[thm]{Lemma}
\newtheorem{prop}[thm]{Proposition}

\newtheorem{cor}[thm]{Corollary}

\newcommand{\prf}{\noindent{\bf Proof.}}
\newcommand{\qed}{\hspace*{\fill}\rule{1.5mm}{2.5mm}}

\title{The $H$-force sets of the graphs satisfying \\ the condition of Ore's theorem}
\author{Xinhong Zhang$^a$
         \thanks{Corresponding author. {\it E-mail address}: xinhongzhang@tyust.edu.cn(X.Zhang).
         Research is supported partially by the National Natural Science Foundation of China (11401353).}, \ \
        Ruijuan Li$^b$\\
\normalsize   $^a${\scriptsize\it Department of Applied Mathematics,
              Taiyuan University of Science and Technology, 030024 Taiyuan, PR China}\\
\normalsize   $^b${\scriptsize\it School of Mathematical Sciences
, Shanxi University, 030006 Taiyuan, PR China}\\}
\begin{document}
\date{}
\maketitle

\begin{abstract}
Let $G$ be a Hamiltonian graph with $n$ vertices. A nonempty vertex set $X\subseteq V(G)$ is called a Hamiltonian cycle enforcing set (in short, an $H$-force set) of $G$ if every $X$-cycle of $G$ (i.e., a cycle of $G$ containing all vertices of $X$) is a Hamiltonian cycle. For the graph $G$, $h(G)$ is the smallest cardinality of an $H$-force set of $G$ and call it the $H$-force number of $G$. Ore's theorem states that the graph $G$ is Hamiltonian if $d(u)+d(v)\geq n$ for every pair of nonadjacent vertices $u,v$ of $G$.  In this paper, we study the $H$-force sets of the graphs satisfying the condition of Ore's theorem, show that the $H$-force number of these graphs is possibly $n$, or $n-2$, or $\frac{n}{2}$ and give a classification of these graphs due to the $H$-force number.
\end{abstract}
\vskip 3mm

\section{Terminology and introduction}

In this paper, we study the simple graphs without loops or no parallel edges. For terminology and notations not defined here we refer the reader to \cite{Bondy}. Let $G$ be a graph with $n$ vertices. The vertex set and the edge set of $G$ are denoted by $V(G)$ and $E(G)$, respectively. For a subset $X$ $\subseteq$ $V(G)$, the cardinality of $X$ is denoted by $|X|$. A subgraph induced by a subset $X$ is denoted by $G[X]$. In addition, $G-X=G[V(G)-X]$.

Let $u$, $v$ be two distinct vertices and $e$ be an edge of $D$. We say that a vertex $v$ is incident to an edge $e$ if $v$ is an endpoint of $e$. If the vertices $u,v$ are incident to the same edge, we say that $u$ and $v$ are {\it adjacent}. Let $S$ be a subset of $V(G)$ or a subgraph of $G$ and $v$ is not in $S$. The set $N_{S}(v)=\{x|\,vx\in E,\,x\in S\}$ and $d_{S}(v)=|N_{S}(v)|$. Let $\delta(G)$ denote the minimum degree of $G$.

A graph $G$ is said to be {\it $k$-connected}, if $V(G)\geq k+1$ and $G-S$ is connected for each $S$ $\subseteq$ $V(G)$ with $|S|\leq k-1$. The subset $S$ $\subseteq$ $V(G)$ is called an {\it independent set of} $G$ if any pair of vertices in $S$ are nonadjacent. In addition, the subset $S$ is called the {\it largest independent set} of $G$ if there is no independent set $S'$ such that $|S'|>|S|$.

For two disjoint graphs $G$ and $H$, $G\vee H$ is the graph that every vertex of $G$ is adjacent to every vertex of $H$. A complete graph with $n$ vertices is denoted by $K_{n}$. $K_{m,n}$ is a complete bipartite graph with two partitions $V_{1},V_{2}$ and $|V_{1}|=m,|V_{2}|=n$.

A cycle of $G$ is called the {\it Hamiltonian cycle}, if it contains all vertices of $G$. The graph $G$ is said to be a {\it Hamiltonian graph} if it has a Hamiltonian cycle. It is well known that the problem to check whether a graph has a Hamiltonian cycle or not is NP-complete. The various kinds of sufficient conditions to imply a graph to be Hamiltonian have been studied widely. More details can be found in [2-5].

Let $G$ be a Hamiltonian graph with $n$ vertices. For a subset $X$ $\subseteq$ $V(G)$, an {\it $X$-cycle} of $G$ is a cycle containing all vertices of $X$.  A nonempty vertex set $X\subseteq V(G)$ is called a {\it Hamiltonian cycle enforcing set} (in short, an {\it $H$-force set}) of $G$ if every $X$-cycle of $G$ is a Hamiltonian cycle. For the graph $G$, $h(G)$ is the smallest cardinality of an $H$-force set of $G$ and call it the {\it $H$-force number} of $G$. The subset $X$ is called the {\it minimum $H$-force set} of $G$ if $X$ is an $H$-force set  with $|X|=h(G)$, i.e., $X$ is an $H$-force set and there is no $H$-force set $X'$ with $|X'|<|X|$.

The $H$-force set and $H$-force number were introduced by Fabrici et al. in \cite{Fabrici}. These concepts play an important role on the Hamiltonian problems. In [6], the authors also investigated the $H$-force numbers of the bipartite graphs, the outerplanar Hamiltonian graphs and so on. Recently, these concepts were generalized to digraphs by Zhang et al. in \cite{X. Zhanglocally} and to hypertournaments by Li et al. in \cite{hypertournament}, and they gave the characterization of the minimum $H$-force sets of locally semicomplete diagraphs and hypertournaments, and obtained their $H$-force numbers.

If $d(u)+d(v)\geq n$ for every pair of nonadjacent vertices $u$ and $v$ of $G$, we say that $G$ satisfies the condition of Ore's theorem. For convenience, we call a graph satisfying the condition of Ore's theorem an OTG. In \cite{Ore}, Ore proved that any OTG is Hamiltonian. In this paper, we study the $H$-force sets and $H$-force number of the OTG's.

The following is an easy observation on the $H$-force sets of a graph.

\begin{prop}\cite{Fabrici} Let $G$ be a Hamiltonian graph. If $C$ is a nonhamiltonian cycle of $G$, then any $H$-force set of $G$ contains a vertex of $V(G)\setminus V(C)$.
\end{prop}

To present Theorem 1.2, we consider a special class of graphs, namely, $$\psi_{2m+1}=\{Z_{m}\vee(K^{c}_{m}+\{u\})\,|\,Z_{m}\; is\; a\; graph\; with\; m\; vertices\}.$$
where $K^{c}_{m}$ is a set of $m$ vertices (also as the complement of the complete graph $K_{m}$) and $u$ is another single vertex, furthermore, the edge set of $Z_{m}\vee(K^{c}_{m}+\{u\})$ consists of $E(Z_{m})$ and $\{xy\,|\,x\in V(Z_{m})\;and\;y\in K^{c}_{m}\cup\{u\}\}$(see Fig.1). In [2], Li et al. investigated the graphs satisfying $d(u)+d(v)\geq n-1$ for every pair of vertices $u,v$ with $d(u,v)=2$ and proved the following.

\begin{thm}\cite{efficient conditions} Let $G$ be a 2-connected graph with $n\geq 3$ vertices. If $d(u)+d(v)\geq n-1$ for every pair of vertices $u,v$ with $d(u,v)=2$, then $G$ is a Hamiltonian graph, unless $n$ is odd and $G\in\psi_{n}$.
\end{thm}

By Theorem 1.2, we obtain immediately the following corollary.

\begin{cor} Let $G$ be a 2-connected graph with $n\geq 3$ vertices. If $d(u)+d(v)\geq n-1$ for every pair of nonadjacent vertices $u,v$ of $G$, then $G$ is a Hamiltonian graph, unless $n$ is odd and $G\in\psi_{n}$.
\end{cor}
\vskip -3mm

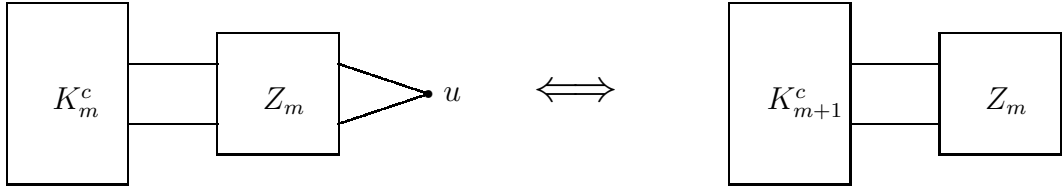
\begin {figure}[h]
\unitlength0.4cm
\begin{center}
\begin{picture}(37,7)
\put(1,0){\framebox(4,6)}
\put(8,1){\framebox(4,4)}
\put(15,3){\circle*{.3}}
\put(15.5,2.75){$u$}
\put(2.5,2.5){$K^{c}_{m}$}
\put(9.5,2.5){$Z_{m}$}
\qbezier(5,4)(6.5,4)(8,4)
\qbezier(5,2)(6.5,2)(8,2)
\qbezier(12,4)(13.5,3.5)(15,3)
\qbezier(12,2)(13.5,2.5)(15,3)
\put(18.5,2.75){\Large{$\Longleftrightarrow$}}

\put(25,0){\framebox(4,6)}
\put(32,1){\framebox(4,4)}
\put(26.25,2.5){$K^{c}_{m+1}$}
\put(33.5,2.5){$Z_{m}$}
\qbezier(29,4)(31.5,4)(32,4)
\qbezier(29,2)(31.5,2)(32,2)
\end{picture}
\caption{$G\in\psi_{2m+1}$}

\end{center} \end{figure}

\section{The H-force sets of the OTG's}

\begin{lem}
Let $G$ be an OTG with $n$ vertices and $X\subseteq V(D)$ be a vertex subset. If $d(u)+d(v)\geq n+1$ for some pair of nonadjacent vertices $u$ and $v$, then

(1) $X$ is an $H$-force set of $G$ if and only if $X$ is an $H$-force set of $G+uv$, and

(2) $X$ is the minimum $H$-force set of $G$ if and only if $X$ is the minimum $H$-force set of $G+uv$, namely, $h(G)=h(G+uv)$.
\end{lem}
{\prf} Firstly, we prove (1). Note that a cycle of the graph $G$ is always a cycle of the graph $G+uv$. It is easy to see that if $X$ is an $H$-force set of $G+uv$, then $X$ is an $H$-force set of $G$.

Suppose that $X$ is an $H$-force set of $G$, but not an $H$-force set of $G+uv$. Let $C$ be the longest nonhamiltonian cycle containing all vertices of $X$ in $G+uv$ and $H$ be the subgraph induced by $V(C)$ in $G$, i.e., $H=G[V(C)]$. Clearly, the cycle $C$ contains the edge $uv$. Assume that there are $a(\geq1)$ vertices apart from $V(C)$. Then $H$ contains $n-a$ vertices. Let $C=v_{1}v_{2}...v_{n-a}v_{1}$ with $u=v_{1}$ and $v=v_{n-a}$. Then $v_{1}v_{2}...v_{n-a}$ is a $(u,v)$-path of $G$ containing all vertices of $X$. Let $S=\{v_{i}\,|\,uv_{i+1}\in E(H)\}$, $T=\{v_{i}\,|\,v_{i}v\in E(H)\}$, where the subscripts are taken modulo $n-a$. Since $v_{n-a}\notin S\cup T$, we have $|S\cup T|<n-a=|V(H)|$.

Recall that $C$ is the longest nonhamiltonian cycle containing all vertices of $X$ in $G+uv$.
For $a>1$, the vertices $u$ and $v$ have no common adjacent vertex apart from $V(C)$. So
$$d_{H}(u)+d_{H}(v)\geq d(u)+d(v)-a\geq n+1-a=n-a+1=|V(H)|+1.$$
For $a=1$, the common adjacent vertex of $u,v$ in $G-C$ is possibly the unique vertex of $G-C$ . So
$$d_{H}(u)+d_{H}(v)\geq d(u)+d(v)-2\geq n+1-2=n-1=|V(H)|.$$
In both cases, $d_{H}(u)+d_{H}(v)\geq |V(H)|.$  Since $$d_{H}(u)+d_{H}(v)=|S|+|T|=|S\cup T|+|S\cap T|\geq |V(H)|$$ and $|S\cup T|<|V(H)|$, then $|S\cap T|>0$. Let $v_{i}\in S\cap T$. Then there is a nonhamiltonian cycle $v_{1}v_{2}...v_{i}v_{n-a}v_{n-a-1}...v_{i+1}v_{1}$ containing all vertices of $X$ in $G$ (see Fig.2).
This contradicts with the fact that $X$ is an $H$-force set of $G$. Thus if the subset $X$ is an $H$-force set of $G$, $X$ is an $H$-force set of $G+uv$.

\begin {figure}[h]
\unitlength0.4cm
\begin{center}
\begin{picture}(27,4)
\put(1,1){\circle*{.3}}
\put(4,1){\circle*{.3}}
\put(7,1){\circle*{.3}}
\put(10,1){\circle*{.3}}
\put(13,1){\circle*{.3}}
\put(16,1){\circle*{.3}}
\put(19,1){\circle*{.3}}
\put(22,1){\circle*{.3}}
\put(25,1){\circle*{.3}}
\put(8,1){$. . .$}
\put(20,1){$. . .$}
\put(0,0){$(u)v_{1}$}
\put(3.75,0){$v_{2}$}
\put(6.5,0){$v_{3}$}
\put(12.75,0){$v_{i}$}
\put(15.25,0){$v_{i+1}$}
\put(21,0){$v_{n-a-1}$}
\put(24.5,0){$v_{n-a}(v)$}
\qbezier(1,1)(4,1)(7,1)
\qbezier(11,1)(15,1)(19,1)
\qbezier(22,1)(23.5,1)(25,1)
\qbezier(1,1)(8,5)(16,1)
\qbezier(13,1)(18,5)(25,1)
\linethickness{1.3pt}
\qbezier(1,1)(4,1)(7,1)
\qbezier(10,1)(11.5,1)(13,1)
\qbezier(16,1)(17.5,1)(19,1)
\qbezier(22,1)(23.5,1)(25,1)
\qbezier(13,1)(18,5)(25,1)
\qbezier(1,1)(8,5)(16,1)
\end{picture}
\caption{A nonhamiltonian cycle $v_{1}v_{2}...v_{i}v_{n-a}v_{n-a-1}...v_{i+1}v_{1}$}

\end{center}
\end{figure}
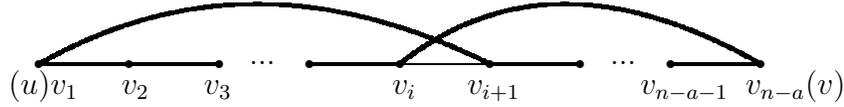

By (1), it is easy to see that (2) holds.                     \qed
\vskip 3mm
Lemma 2.1 motivates the following definition. The {\it weak closure} of an OTG $G$ is the graph obtained from $G$ by recursively joining pairs of nonadjacent vertices whose degree sum is at least $n+1$ until no such pair remains. We denote the  weak closure of $G$ by $C_{w}(G)$. The idea of weak closure is followed by Bondy and Chv\'{a}tal's introduction of the {\it closure} of a graph. If we join pairs of nonadjacent vertices with the degree sum at least $n$ until no such pair remains, then we obtain the closure of a graph, denoted by $C(G)$. Bondy and Chv\'{a}tal proved that $C(G)$ is well defined. That means, if $G_{1}$ and $G_{2}$ are two graphs obtained from $G$ by recursively joining pairs of nonadjacent vertices whose degree sum is at least $n$ until no such pair remains, then $G_{1}=G_{2}$. By the similar argument, we can see that $C_{w}(G)$ is well defined.

Lemma 2.1 implies that the $H$-force set and $H$-force number of an OTG $G$ are the $H$-force set and $H$-force number of the $C_{w}(G)$, respectively. We can obtain the $H$-force set and $H$-force number of $G$ by studying the weak closure $C_{w}(G)$ of $G$. Clearly, for an OTG $G$,  the weak closure $C_{w}(G)$ is a complete graph or a graph satisfying $d_{C_{w}(G)}(u)+d_{C_{w}(G)}(v)=n$ for every pair of nonadjacent vertices $u$ and $v$.

\begin{thm}
Let $G$ be an OTG with $n\geq5$ vertices and $X$ be the minimum $H$-force set of $G$. Let $C_{w}(G)$ be the {\it weak closure} of $G$ and $S$ be the largest independent set of $C_{w}(G)$. Then

(1) the $H$-force number $h(G)=n-2$, $\frac{n}{2}$, or $n$, and

(2) $$
X=\left\{
\begin{array}{ll}
$V(G)-\{x,y\}$,&\mbox {if $h(G)=n-2$,}\\
$S$,&\mbox {if $h(G)=\frac{n}{2}$,}\\
$V(G)$,&\mbox {if $h(G)=n$,}\\
\end{array}
\right.
$$

\noindent where $x,y\in V(G)$ with $d_{C_{w}(G)}(x)=n-1$ and $d_{C_{w}(G)}(y)=n-1$.
\end{thm}

{\prf} According to Lemma 2.1, it is sufficient to consider the minimum $H$-force set and $H$-force number of the weak closure $C_{w}(G)$. For the convenience, let $G_{w}=C_{w}(G)$.

Obviously, the minimum $H$-force set of $G$ is $V(G)$ and $h(G)=n$ if $G_{w}$ is a complete graph .

So assume that $G_{w}$ is not a complete graph and satisfies $d_{G_{w}}(u)+d_{G_{w}}(v)=n$ for every pair of nonadjacent vertices $u$ and $v$ in $G_{w}$. Clearly, $2\leq\delta(G_{w})\leq\frac{n}{2}$.

We consider the following two cases.

\vskip 3mm

{\it Case 1:} $2\leq\delta(G_{w})<\frac{n}{2}$.

\vskip 3mm
Let $\delta(G_{w})=a$ and $d_{G_{w}}(u)=\delta(G_{w})=a$ for some $u\in V(G_{w})$. Then $d_{G_{w}}(u)<\frac{n}{2}$. We use the following notations: $$U=\{v\in V-u:\,uv\in E(G_{w})\}, \;W=V-U-\{u\}.$$ Obviously, $|U|=a$, $|W|=n-a-1$ (see Fig. 3).

\begin {figure}[h]
\unitlength0.4cm \begin{center}
\begin{picture}(33,7)
\put(3,1){\framebox(3,5)}
\put(9,1){\framebox(3,5)}
\put(1,3){\circle*{.3}}
\put(4.5,5.5){\circle*{.3}}
\put(4.5,4.5){\circle*{.3}}
\put(4.5,2.5){\circle*{.3}}
\put(4.5,1.5){\circle*{.3}}
\put(10.5,5){\circle*{.3}}
\put(0,2.75){$u$}
\put(3.75,5.25){$x$}
\put(3.75,4.25){$y$}
\put(11,4.75){$w$}
\put(3.75,2.25){$q$}
\put(3.75,1.25){$p$}
\put(4.25,0){$U$}
\put(10.25,0){$W$}

\qbezier(1,3)(2,3.75)(3,4.5)
\qbezier(1,3)(2,2.25)(3,1.5)
\qbezier(4.5,5.5)(7.5,5.25)(10.5,5)
\qbezier(4.5,4.5)(7.5,4.75)(10.5,5)

\put(18,1){\framebox(3,5)}
\put(25,1){\framebox(3,5)}
\put(23,3){\circle*{.3}}
\put(23,5){\circle*{.3}}
\put(22.6,2.2){$y$}
\put(22.6,5.5){$x$}
\put(19.25,0){$U'$}
\put(26.25,0){$W$}
\put(18.5,3){\footnotesize{$K_{a-1}$}}
\put(25,3){\footnotesize{$K_{n-a-1}$}}
\qbezier(21,4.5)(22,4.75)(23,5)
\qbezier(21,3.5)(22,3.25)(23,3)
\qbezier(25,4.5)(24,4.75)(23,5)
\qbezier(25,3.5)(24,3.25)(23,3)
\qbezier(23,3)(23,4)(23,5)
\put(7,-1){$G_{w}$}
\put(23,-1){$G_{11}$}

\end{picture}
\caption{$G_{w}$ and $G_{11}\in\varphi_{1}$}

\end{center}
\end{figure}
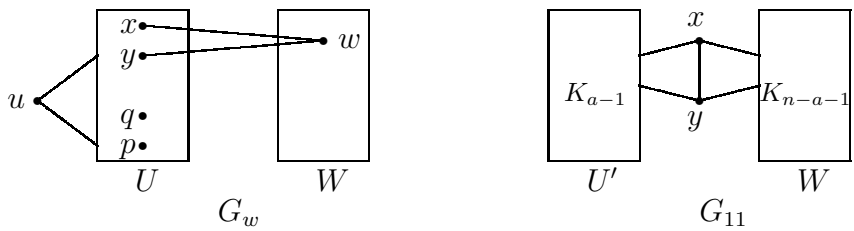

\vskip 3mm

\noindent {\bf Claim A:}    (i) For each vertex $w\in W$, $d_{G_{w}}(w)=n-a$.

(ii) $W$ induces a complete subgraph.

(iii) For each vertex $w\in W$, $d_{U}(w)=2$.

Let $w\in W$, and $x,y\in U$ be the vertices adjacent to $w$ in $U$.

(iv) If $U-\{x,y\}\neq \phi$, then $d_{G_{w}}(p)=a$ for each vertex $p\in U-\{x,y\}$.

(v) If $U-\{x,y\}\neq \phi$, then $U-\{x,y\}$ induces a complete subgraph.

(vi) $x$ and $y$ are adjacent to every vertex of $U-\{x,y\}$ if $U-\{x,y\}\neq\phi$.

(vii) For any vertex $q\in U-\{x,y\}$, $q$ has no neighbour in $W$.

(viii) $x$ and $y$ are adjacent to every vertex of $W$.

(ix) $d_{G_{w}}(x)=d_{G_{w}}(y)=n-1$ and $U$ induces a complete subgraph.
\vskip 3mm
\noindent{\it Proof.}   (i) For the vertex $w\in W$, since $w$ is not adjacent to $u$, we see that $d_{G_{w}}(w)=n-d_{G_{w}}(u)=n-a$.

(ii) For the distinct vertices $w,w'\in W$, we have $d_{G_{w}}(w)=d_{G_{w}}(w')=n-a>\frac{n}{2}$, and hence $d_{G_{w}}(w)+d_{G_{w}}(w')>n$. So $w$ and $w'$ are adjacent. By the choice of $w,w'$, $W$ induces a complete subgraph.

(iii) Since $|W|=n-a-1$ and $W$ is a complete subgraph, we have $$d_{U}(w)=d_{G_{w}}(w)-(n-a-2)=(n-a)-(n-a-2)=2.$$

(iv) Since the vertex $p$ is not adjacent to $w$ and $d_{G_{w}}(w)=n-a$, we have $d_{G_{w}}(p)=n-d_{G_{w}}(w)=n-(n-a)=a$.

(v) If $p,q$ are a pair of nonadjacent vertices in $U-\{x,y\}$, then $d_{G_{w}}(p)+d_{G_{w}}(q)=n$. But $d_{G_{w}}(p)+d_{G_{w}}(q)=2a<n$. So $p$ and $q$ are adjacent and $U-\{x,y\}$ induces a complete subgraph.

(vi) Suppose not. Let $p$ be a vertex of $U-\{x,y\}$ such that $x$ is not adjacent to $p$. By (iv), we have $d_{G_{w}}(x)=n-d_{G_{w}}(p)=n-a$. For a vertex $w\in W$, $d_{G_{w}}(w)=n-a$ by (i). Then $d_{G_{w}}(x)+d_{G_{w}}(w)>n$. So $x$ is adjacent to every vertex of $W$.

Recall that $|W|=n-a-1$ and $d_{G_{w}}(x)=n-a$. Hence $x$ has no neighbour in $U-\{x\}$. In particular, $x$ is not adjacent to $y$. So $d_{G_{w}}(y)=n-d_{G_{w}}(x)=a$.  Since the degree of any vertex of $U-\{x,y\}$ is also $a$, the vertex $y$ is adjacent to every vertex of $U-\{x,y\}$. So, $w$ is the only adjacent vertex of $y$ in $W$.  Therefore, $$N_{G_{w}}(x)=\{u\}\cup W,\quad N_{G_{w}}(y)=\{u\}\cup (U-\{x,y\})\cup \{w\}.$$ In other words, $\{x\}\cup (W-\{w'\})\subseteq N_{G_{w}}(w')$ for any vertex $w'\in W-\{w\}$. By $|\{x\}\cup(W-\{w'\})|=n-a-1$ and $d_{G_{w}}(w')=n-a$, we see that every vertex of $W-\{w\}$ has exactly an adjacent vertex in $U-\{x,y\}$.

Let $q\in U-\{x,y\}$ be arbitrary. By the arguments above, we see that $\{u\}\cup(U-\{x,q\})\subseteq N_{G_{w}}(q)$. Since $d_{G_{w}}(q)=a$ (see (iv)) and$|\{u\}\cup(U-\{x,q\})|=a-1$,
we have that $q$ has an adjacent vertex in $W-\{w\}$. This means every vertex of $U-\{x,y\}$ has exactly an adjacent vertex in $W-\{w\}$.

So, $|U-\{x,y\}|=|W-\{w\}|$. But $a-2<n-a-2$, a contradiction. Thus $x$ is adjacent to every vertex of $U-\{x,y\}$. Similarly, $y$ is adjacent to every vertex of $U-\{x,y\}$.

(vii) For $q\in U-\{x,y\}$, the claims (iv), (v) and (vi) imply that $q$ is not adjacent to any vertex of $W$.

(viii) The claims (iii) and (vii) show that every vertex of $W$ has two adjacent vertices in $U$, which are just $x$ and $y$.

(ix) By (vi) and (viii), $d_{G_{w}}(x)+d_{G_{w}}(y)\ge 2(|U-\{x,y\}|+|W|+|\{u\}|)=2(n-2)>n$. Then $x$ and $y$ are adjacent and hence (ix) holds.       \qedc

\vskip 5mm

Now, $G_{w}$ is isomorphic to $G_{11}=K_2\vee (K_{a-1}+K_{n-a-1})$ (see Fig.3), where $U'=(U-\{x,y\})\cup\{u\}$. In particular, $d_{G_{w}}(x)=d_{G_{w}}(y)=n-1$, $U'$ and $W$ induce complete subgraphs with $|U'|=a-1$ and $|W|=n-a-1$.

We shall show that $X=V-\{x,y\}=U'\cup W$ is the minimum $H$-force set and hence $h(G_{w})=n-2$.

Obviously, $X$ is an $H$-force set because a cycle containing all vertices of $X$ must encounter $x$ and $y$. Suppose to the contrary that $X'$ is an $H$-force set of $G_{w}$ with $|X'|<|X|$. Clearly, there exists a vertex $z\in U'$ or $z\in W$  such that $z\notin X'$. Without loss of generality, assume $z\in U'$. Consider the subgraph $G_{w}-z$.
For any pair of nonadjacent vertices $u\in U'-z,v\in W$ of $G_{w}-z$, we have $$d_{G_{w}-z}(u)+d_{G_{w}-z}(v)=d_{G_{w}}(u)-1+d_{G_{w}}(v)=n-1.$$
Then $G_{w}-z$ is an OTG and hence $G_{w}-z$ is a Hamiltonian graph. In other words, there is a nonhamiltonian cycle containing all vertices of $X'$ in $G_w$, a contradiction. Thus there doesn't exist such an $H$-force set $X'$ with $|X'|<|X|$. So $X$ is the minimum $H$-force set and $h(G_{w})=h(G)=n-2$.

\vskip 3mm {\it Case 2:}  $\delta(G_{w})=\frac{n}{2}$ and $n$ is even.
\vskip 3mm
{\it Subcase 2.1:}
There is $v\in V(G_{w})$ such that $G_{w}-v$ is not a Hamiltonian graph.
\vskip 3mm

Let $H=G_{w}-v$. For arbitrary $x,y\in V(H)$, we have $$d_{H}(x)+d_{H}(y)\geq d_{G_{w}}(x)-1+d_{G_{w}}(y)-1=n-2=|V(H)|-1.$$
 By Corollary 1.3, either $H$ is not a 2-connected graph or $H\in \psi_{n-1}$(see Fig.1).

\vskip 3mm

{\it Subcase 2.1.1:} $H=G_{w}-v$ is not a 2-connected graph.

\vskip 3mm

Then there exists a vertex $u$ such that $H-u$ is not connected.

\vskip 3mm

\noindent {\bf Claim B:} (i) For any vertex $x\in V(G_{w})$, if $d_{G_{w}}(x)>\frac{n}{2}$, then $d_{G_{w}}(x)=n-1$.

(ii) $H-u$ has exactly two components, say $U$ and $W$.

(iii) For any vertices $x\in U$, $y\in W$, $d_{G_{w}}(x)=d_{G_{w}}(y)=\frac{n}{2}$.

(iv) $|U|=|W|=\frac{n}{2}-1$.

(v) Both $U$ and $W$ induce complete subgraphs, and $d_{G_{w}}(u)=d_{G_{w}}(v)=n-1$.

\vskip 3mm

\noindent{\it Proof.} (i)  Recall that $\delta(G_{w})=\frac{n}{2}$, and $d_{G_{w}}(u)+d_{G_{w}}(v)=n$ for nonadjacent vertices $u,v$. It is clear that (i) holds.

(ii) Suppose to the contrary that there are three components $U,W$ and $U'$ in $H-u$. So $|U|+|W|\leq n-3$. For the vertices $x\in U$, $y\in W$, we see that
$$n=d_{G_{w}}(x)+d_{G_{w}}(y)=d_{U}(x)+d_{W}(y)+d_{\{u,v\}}(x)+d_{\{u,v\}}(y)$$   $$\qquad\qquad\quad\leq |U|-1+|W|-1+4=|U|+|W|+2\leq n-3+2=n-1,$$
a contradiction. Thus $H-u$ has two components.

(iii) For $x\in U,\,y\in W$, we have $d_{G_{w}}(x)\neq n-1$ and $d_{G_{w}}(y)\neq n-1$. By (i), $d_{G_{w}}(x)=d_{G_{w}}(y)=\frac{n}{2}$.

(iv) Note that $|U\cup W|=n-2$. Assume without loss of generality that $|U|\leq\frac{n}{2}-1$. For any vertex $x\in U$, we have $$\frac{n}{2}=d_{G_{w}}(x)\leq|U|-1+2=|U|+1.$$
 Then $|U|\geq \frac{n}{2}-1$. So $|U|=|W|=\frac{n}{2}-1$.

(v) For any vertex $x\in U$, $d_{G_{w}}(x)=\frac{n}{2}\leq d_{U}(x)+2\le |U|-1+2=\frac{n}{2}-1+1=\frac{n}{2}$. Then $d_{G_{w}}(x)=d_{U}(x)+2=(|U|-1)+2$.  This implies that $U$ induces a complete subgraph and $x$ is adjacent to $u,v$. Similarly, $W$ induces a complete subgraph and every vertex of $W$ is adjacent to $u,v$.

Furthermore, all vertices of $U$ and $W$ are the neighbours of $u$ and $v$. So $d_{G_{w}}(u)\geq n-2$, $d_{G_{w}}(v)\geq n-2$ and hence $d_{G_{w}}(u)+d_{G_{w}}(v)>n$. It means that $u$ and $v$ are adjacent. Thus $d_{G_{w}}(u)=d_{G_{w}}(v)=n-1$.           \qedc

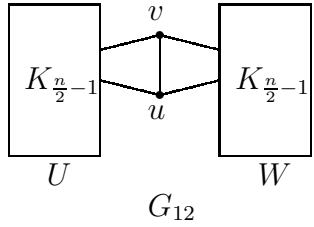
\begin {figure}[h]
\unitlength0.4cm
\begin{center}
\begin{picture}(13,7)

\put(1,1){\framebox(3,5)}
\put(8,1){\framebox(3,5)}
\put(6,3){\circle*{.3}}
\put(6,5){\circle*{.3}}
\put(1.5,3.2){$K_{\frac{n}{2}-1}$}
\put(8.5,3.2){$K_{\frac{n}{2}-1}$}
\put(5.6,2.2){$u$}
\put(5.6,5.5){$v$}
\put(2.25,0){$U$}
\put(9.25,0){$W$}
\put(5.6,-1){$G_{12}$}
\qbezier(4,4.5)(5,4.75)(6,5)
\qbezier(4,3.5)(5,3.25)(6,3)
\qbezier(8,4.5)(7,4.75)(6,5)
\qbezier(8,3.5)(7,3.25)(6,3)
\qbezier(6,3)(6,4)(6,5)
\end{picture}
\caption{$G_{12}\in\varphi_{1}$}
\end{center} \end{figure}

\vskip 3mm
In subcase 2.1.1, $G_{w}\cong G_{12}$ (see Fig.4), in which $U$ and $W$ are complete subgraphs with $|U|=|W|=\frac{n}{2}-1$. By the same argument of the case when $2\leq\delta(G_{w})<\frac{n}{2}$ and $G_{w}\cong G_{11}$, we also have $V-\{u,v\}$ is the minimum $H$-force set and $h(G_{w})=h(G)=n-2$.

\vskip 3mm
{\it Subcase 2.1.2:} $H=G_{w}-v\in \psi_{n-1}$.
\vskip 3mm

Let $H=G_{w}-v=Z_{m}\vee(K_{m+1}^{c})$, where $Z_{m}$ is a graph with $m$ vertices. Clearly, $|V(H)|=n-1=2m+1$, namely, $n=2m+2$. Recall that $n\geq5$. So $m\geq2$.

\vskip 3mm
\noindent {\bf Claim C:}  If there is at least one edge in $Z_{m}$, then

(i) $Z_{m}$ is a complete subgraph.

(ii) The vertex $v$ is adjacent to all vertices of $Z_{m}$.

(iii) The vertex $v$ is adjacent to all vertices of $K_{m+1}^{c}$.

(iv) $G_{w}\cong G_{21}=K_{m+1}^c\vee K_{m+1}$ (see Fig.5). Furthermore, $X=V(K_{m+1}^c)$  is the minimum $H$-force set and $h(G)=h(G_w)=\frac{n}{2}$.
\vskip 3mm

\noindent{\it Proof.} (i) Let $x$ be an endpoint of the edge of $Z_{m}$. Then $d_{G_{w}}(x)\geq m+2$. For arbitrary $y\in Z_{m}-x$, since $d_{G_{w}}(y)\geq m+1$, we have $d_{G_{w}}(x)+d_{G_{w}}(y)\geq 2m+3>n$. This implies that $x$ is adjacent to the remaining vertices of $Z_{m}$. Furthermore, for any pair of $p,q\in Z_{m}-x$, we have $d_{G_{w}}(p)\geq m+2, d_{G_{w}}(q)\geq m+2$, and hence $d_{G_{w}}(p)+d_{G_{w}}(q)\geq 2m+4>n$, which implies that $p$ and $q$ are adjacent. Thus $Z_{m}$ is a complete subgraph.

(ii) For a vertex $z\in Z_{m}$, $$d_{G_{w}}(v)+d_{G_{w}}(z)\geq \frac{n}{2}+(n-2)=m+1+2m>2m+2=n.$$
Thus $v$ is adjacent to all vertices of $Z_{m}$.

(iii) Since $\delta(G_{w})=\frac{n}{2}$, $d_{G_{w}}(x)\geq \frac{n}{2}=m+1=|V(Z_m)|+1$ for any vertex $x\in K_{m+1}^{c}$. So $x$ is adjacent to $v$. Thus  $v$ is adjacent to all vertices of $K_{m+1}^{c}$.

\vskip 3mm

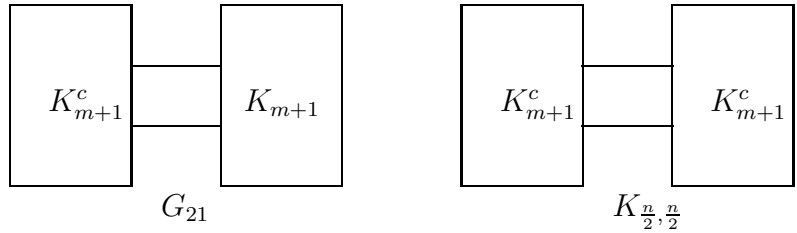
\begin {figure}[h]
\unitlength0.4cm
\begin{center}
\begin{picture}(29,7)
\put(1,1){\framebox(4,6)}
\put(8,1){\framebox(4,6)}
\put(16,1){\framebox(4,6)}
\put(23,1){\framebox(4,6)}
\put(2.25,3.5){$K^{c}_{m+1}$}
\put(17.25,3.5){$K^{c}_{m+1}$}
\put(8.75,3.5){$K_{m+1}$}
\put(24.25,3.5){$K^{c}_{m+1}$}
\qbezier(5,5)(6.5,5)(8,5)
\qbezier(20,5)(21.5,5)(23,5)
\qbezier(5,3)(6.5,3)(8,3)
\qbezier(20,3)(21.5,3)(23,3)
\put(6,0){$G_{21}$}
\put(21,0){$K_{\frac{n}{2},\frac{n}{2}}$}
\end{picture}
\caption{$G_{21}\in\varphi_{2}$ and $K_{\frac{n}{2},\frac{n}{2}}\in\varphi_{2}$}
\end{center}
\end{figure}

(iv) By the arguments above, it is obvious $G_{w}\cong G_{21}=K_{m+1}^c\vee K_{m+1}$.

Let $X=V(K^{c}_{m+1})$. For any $x\in X$, it is obvious that $G_{w}-x$ is an OTG and hence it is Hamiltonian. By Proposition 1.1, the minimum $H$-force set contains all vertices of $V(K^{c}_{m+1})$. In addition, it is clear that $X=V(K^{c}_{m+1})$ is an $H$-force set of $G_{w}$. So $X=V(K^{c}_{m+1})$ is the minimum $H$-force set and $h(G)=h(G_{w})=\frac{n}{2}$.
\qedc

\vskip 3mm

\noindent {\bf Claim D:}  If there is no edge in $Z_{m}$, namely, $Z_{m}=K^{c}_{m}$, then

(i) $v$ has no neighbour in $Z_{m}$.

(ii) $v$ is adjacent to all vertices of $K_{m+1}^{c}$.

(iii) $G_{w}\cong K_{\frac{n}{2},\frac{n}{2}}$(see Fig.5).  Furthermore, $X=V(K^{c}_{m+1})$ is the minimum $H$-force set and $h(G)=h(G_w)=\frac{n}{2}$.
\vskip 3mm

\noindent{\it Proof.} (i) Suppose not. Let $v$ is adjacent to $x$ for some $x\in V(Z_{m})$. Let $y$ be another vertex of $Z_{m}$. Then $d_{G_{w}}(x)+d_{G_{w}}(y)\geq m+2+m+1>n$ and hence $x$ and $y$ are adjacent, which contradicts with $Z_{m}=K_{m}^{c}$.

(ii)  Since $\delta(G_{w})=\frac{n}{2}$, we have $v$ is adjacent to all vertices of $K_{m+1}^{c}$.

(iii) By (i) and (ii), it is easy to see that $G_{w}\cong K_{\frac{n}{2},\frac{n}{2}}$.  Clearly, the largest independent set $V(K^{c}_{m+1})$ is the minimum $H$-force set and $h(G)=h(G_{w})=\frac{n}{2}$.      \qedc

\vskip 3mm
{\it Subcase 2.2:}  For any vertex $v\in V(G_{w})$,  the subgraph $G_{w}-v$ is Hamiltonian.

 In this subcase, Proposition 1.1 yields that the minimum $H$-force set is $V(G)$ and $h(G)=h(G_{w})=n$.
\qed

\vskip 3mm

In \cite{Dirac}, Dirac proved that Theorem 2.3.

\begin{thm}\cite{Dirac}
If the minimum degree $\delta(G)$ of $G$ is no less than $\frac{n}{2}$, then $G$ is Hamiltonian.
\end{thm}

Clearly, the graph satisfying the condition of Dirac's theorem is an OTG. By Theorem 2.2, we obtain immediately the following results on the $H$-force sets of graphs satisfying the condition of Dirac's theorem.

\begin{cor}
Let $G$ be a graph satisfying $\delta(G)\geq\frac{n}{2}$
and $X$ be the minimum $H$-force set of $G$. Let $G_{w}=C_{w}(G)$ be the weak closure of $G$ and $S$ be the largest independent set of $G_{w}$. Then

(1) The $H$-force number $h(G)=n-2$, $\frac{n}{2}$, or $n$.

(2) $$
X=\left\{
\begin{array}{ll}
$V(G)-\{x,y\}$,&\mbox {if $h(G)=n-2$,}\\
$S$,&\mbox {if $h(G)=\frac{n}{2}$,}\\
$V(G)$,&\mbox {if $h(G)=n$,}\\
\end{array}
\right.
$$

\noindent where $x,y\in V(G)$ with $d_{G_{w}}(x)=n-1$ and $d_{G_{w}}(y)=n-1$.
\end{cor}
\section{The classification of the OTG's}

To present Theorem 3.1, we define several special classes of graphs, namely,
\begin{eqnarray*}
\varphi_{1}&=&\{K_2\vee (K_{m}+K_{n-m-2})\,|\,1\leq m< \frac{n}{2}\}\, (\mbox{see Fig.3 and Fig.4}),\\
\varphi_{2}&=&\{G_{21},\,K_{\frac{n}{2},\frac{n}{2}}\}\,(\mbox{see Fig.5}),\\
\varphi_{3}&=&\{\{Z_{n-m}\vee K_{m}\,|\, n \mbox{ is even}, 0\leq m< \frac{n}{2}\mbox{ and } d(x)=\frac{n}{2}\mbox{ for }x\in Z_{n-m} \}-\{G_{12},\;K_{\frac{n}{2},\frac{n}{2}}\}\}\\
&&\cup \,\{G\,|\, G \mbox{ is a complete graph}\}\,(\mbox{see\,Fig.6}),
\end{eqnarray*}

\noindent where $Z_{n-m}$ is a graph with $n-m$ vertices.
\vskip 3mm

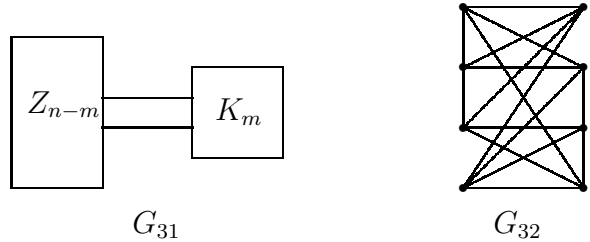
\begin {figure}[h] \unitlength0.4cm
\begin{center}
\begin{picture}(22,9)
\put(1,1){\framebox(3,5)}
\put(7,2){\framebox(3,3)}
\put(16,1){\circle*{.3}}
\put(16,3){\circle*{.3}}
\put(16,5){\circle*{.3}}
\put(16,7){\circle*{.3}}
\put(20,1){\circle*{.3}}
\put(20,3){\circle*{.3}}
\put(20,5){\circle*{.3}}
\put(20,7){\circle*{.3}}
\put(1.5,3.5){$Z_{n-m}$}
\put(7.75,3.25){$K_{m}$}
\qbezier(4,4)(5.5,4)(7,4)
\qbezier(4,3)(5.5,3)(7,3)
\qbezier(16,3)(16,5)(16,7)
\qbezier(20,1)(20,3)(20,5)
\qbezier(16,7)(18,7)(20,7)
\qbezier(16,5)(18,5)(20,5)
\qbezier(16,3)(18,3)(20,3)
\qbezier(16,1)(18,1)(20,1)
\qbezier(16,7)(18,6)(20,5)
\qbezier(16,7)(18,4)(20,1)
\qbezier(16,5)(18,6)(20,7)
\qbezier(16,3)(18,5)(20,7)
\qbezier(16,3)(18,2)(20,1)
\qbezier(16,1)(18,4)(20,7)
\qbezier(16,1)(18,3)(20,5)
\qbezier(16,1)(18,2)(20,3)
\put(5,-0.5){$G_{31}$}
\put(17,-0.5){$G_{32}$}
\end{picture}
\caption{$G_{31}\in\varphi_{3}$ and $G_{32}\in\varphi_{3}$}
\end{center}
\end{figure}

\vskip 3mm

Let $\mathcal{C}_{w}$ be the class of the weak closures of all OTG's. Clearly, $\varphi_{1},\varphi_{2},\varphi_{3}$ are all the subsets of $\mathcal{C}_{w}$. By Theorem 2.2, we can give a partition of the class $\mathcal{C}_{w}$ and obtain the following theorem.

\begin{thm} The class $\mathcal{C}_{w}$ of the weak closures of all OTG's has a partition
$$\mathcal{C}_{w}=\varphi_{1}\,\dot{\cup}\,\varphi_{2}\,\dot{\cup}\,\varphi_{3}$$
due to the $H$-force set, where $\varphi_{1},\varphi_{2},\varphi_{3}$ are defined above. Let $G_{w}\in \mathcal{C}_{w}$ be arbitrary. Then

(1) $h(G_{w})=n-2$ if and only if $G_{w}\in\varphi_{1}$.

(2) $h(G_{w})=\frac{n}{2}$ if and only if $G_{w}\in\varphi_{2}$.

(3) $h(G_{w})=n$ if and only if $G_{w}\in\varphi_{3}$.

\end{thm}

{\prf} By the case 1 and subcase 2.1.1 of Theorem 2.2, the statement (1) holds. Also by the subcase 2.1.2, the statement (2) holds. We show (3) as follows.

``$\Longrightarrow$" By the proof of Theorem 2.2, we know that $h(G_{w})=n-2$ if $2\leq\delta(G_{w})<\frac{n}{2}$. So if $h(G_{w})=n$, then $\delta(G_{w})=\frac{n}{2}$ ($n$ is even) or $G_{w}$ is a complete graph. In the latter case, $G_{w}\in \varphi_{3}$ and we are done.

Assume that $G_{w}$ is not a complete graph. Then $\delta(G_{w})=\frac{n}{2}$ and $d_{G_{w}}(u)+d_{G_{w}}(v)=n$ for every pair of nonadjacent vertices $u,v$.  So the degree of any vertex $u\in V(G_w)$ is either $\frac{n}{2}$ or $n-1$ . Set $U=\{v|\:d_{G_{w}}(v)=n-1\}$ and $m=|U|$.

When $U\ne \emptyset$, we claim that $1\leq m<\frac{n}{2}$. Clearly, $1\le m\le \frac{n}{2}$. Note that $d_{G_{w}}(w)=\frac{n}{2}$ for any vertex $w\in V(G_{w})-U$. If $m=\frac{n}{2}$, then $G_{w}\cong G_{21}$ and $h(G_{w})=\frac{n}{2}$, a contradiction.  So $1\leq m<\frac{n}{2}$.   Since $h(G_{12})=n-2$, we have $G_{w}\ncong G_{12}$. Thus, $G_{w}\in\varphi_{3}$.

When $U$ is an empty set, then $m=0$ and $G_{w}$ is a $\frac{n}{2}$-regular graph. Since $h(K_{\frac{n}{2},\frac{n}{2}})=\frac{n}{2}$, we have $G_{w}\ncong K_{\frac{n}{2},\frac{n}{2}}$. Thus, $G_{w}\in\varphi_{3}$.

``$\Longleftarrow$" It is not difficult to check that $\varphi_{3}\cap\varphi_{1}=\emptyset$ and $\varphi_{3}\cap\varphi_{2}=\emptyset$. Let $G_{w}\in\varphi_{3}$ be a weak closure of an OTG. So $h(G_{w})\neq n-2 \;or \;\frac{n}{2}$. According to Theorem 2.2, $h(G_{w})=n$.
\qed

\vskip 3mm
To present Corollary 3.2, we define several special classes of graphs, namely,
$$\vartheta_{i}=\{G\mbox{ is an OTG}, G_{w}\in \varphi_{i}\}, \quad i=1,2,3.$$

Let $\vartheta$ be the class of all OTG's. Clearly, $\vartheta_{1},\vartheta_{2},\vartheta_{3}$ are all the subsets of $\vartheta$. By Theorem 3.1, we can give a partition of the class $\vartheta$ and obtain the following Corollary.

\begin{cor} The class $\vartheta$ of all OTG's has a partition $$\vartheta=\vartheta_{1}\,\dot{\cup}\,\vartheta_{2}\,\dot{\cup}\,\vartheta_{3}$$
due to the weak closure, where $\vartheta_{1},\vartheta_{2},\vartheta_{3}$ are defined above. Let $G\in \vartheta$ be arbitrary. Then

$$
h(G)=\left\{
\begin{array}{ll}
n-2,&\mbox { if  } G\in \vartheta_{1},\\
\frac{n}{2},&\mbox { if  } G\in \vartheta_{2},\\
n,&\mbox { if  } G\in \vartheta_{3}.\\
\end{array}
\right.
$$

\end{cor}
\section{The Algorithm to check the minimum $H$-force set of an OTG}

\noindent {\bf Input:} A graph $G=(V,E)$

\noindent {\bf Output:} If $G$ is an OTG, return the $H$-force number $h(G)$ and the minimum $H$-force set $X$. If $G$ is not an OTG, return ``not an OTG".

1. For each pair of nonadjacent vertices $u, v$, check $d(u)+d(v)$. If there is $uv\notin E$ with $d(u)+d(v)<n$, then return ``not an OTG".

2. Make the weak closure of $G$, namely, $G_{w}$.

3. For any $x\in V(G)$, check $d_{G_{w}}(x)$.

(1) If there are exactly two vertices $u$ and $v$ with $d_{G_{w}}(u)=d_{G_{w}}(v)=n-1$ and $G_{w}\in \varphi_{1}$, return $h(G)=n-2$ and $X=V(G)-\{u,v\}$.

If there are exactly two vertices $u$ and $v$ with $d_{G_{w}}(u)=d_{G_{w}}(v)=n-1$ and $G_{w}\in \varphi_{3}$, return $h(G)=n$ and $X=V(G)$.

(2) If there are exactly $\frac{n}{2}$ vertices with degree $n-1$, return $h(G)=\frac{n}{2}$. Set $X=\{x|\;d(x)<n-1\}$. Return $X$.

(3) If there are $a$ vertices with degree $n-1$ where $2<a<\frac{n}{2}$ or $a=1$ or $a=n$, return $h(G)=n$ and $X=V(G)$.

(4) If any $v\in V(G_{w})$, $d_{G_{w}}(v)=\frac{n}{2}$ and $G_{w}\cong K_{\frac{n}{2},\frac{n}{2}}$, return $h(G)=\frac{n}{2}$. Find a partition set $X$ of $K_{\frac{n}{2},\frac{n}{2}}$, return $X$.

(5) If any $v\in V(G_{w})$, $d_{G_{w}}(v)=\frac{n}{2}$ and $G_{w}\ncong K_{\frac{n}{2},\frac{n}{2}}$, return $h(G)=n$ and $X=V(G)$.

\vskip 3mm

\begin{thm}
For an OTG, the minimum $H$-force set and the $H$-force number can be found in time O($n^{3}$).

\end{thm}

{\prf} The complexity follows from the fact that step 1 can be performed in time $O(n^{2})$, step 2 in time $O(n^{3})$ and step 3 in time $O(n^2)$.

\qed

\end{document}